\documentstyle[12pt,amssymb]{article}
\begin{document}
\newtheorem{theorem}{Theorem}
\newtheorem{definition}[theorem]{Definition}
\newtheorem{corollary}[theorem]{Corollary}
\newcommand{\bproof}{\noindent{\bf Proof: }}
\newcommand{\eproof}{\hfill $\Box$\par}
\newtheorem{remark}{Remark}
\renewcommand{\theequation}{\thesection.\arabic{equation}}
\setcounter{section}{0}

\title{A note on unconditional structures\\ in weak Hilbert spaces.}
\author{N. J. Nielsen}
\date{ }
\maketitle
\vspace{1cm}

\begin{abstract}
We prove that if a non-atomic separable Banach lattice
in a weak Hilbert space, then it is lattice isomorphic to $L_2(0,1)$.
\end{abstract}

\section*{Introduction}

This note is to be considered as an addendum to \cite{NTJ}, where an
extensive study of Banach lattices, which are weak Hilbert spaces, was
made. We prove that if a non-atomic separable Banach lattice is a weak
Hilbert space then it is lattice isomorphic to $L_2(0,1)$. The result is
a consequence of Theorem 3.11 in [3] together with an easy, short
argument.

We believe that it may have some impact om the study of unconditional
structures in weak Hilbert spaces.

For the convenience of the reader we have in section 1 given the
definition of a weak Hilbert space and formulated the special case of
\cite{NTJ}, Theorem 3.11, which is needed to prove our result.
\vspace{.5cm}

\section{Notation and terminology}

In this note we shall use the notation and terminology commonly used in
Banach space theory as it appears in \cite{LT1} and \cite{LT2}.

If $E$ is a finite dimensional Banach space we denote the Banach-Mazur
distance between $E$ and the Hilbert space by\ d(E),\  and if\ $(x_j)$\ is a
sequence in a Banach space\ $X$\ we let\ $[x_j]$\ denote the closed linear
span of\ $(x_j)$.\  Given a finite set\ $A$\ we let\ $|A|$\ denote the
cardinality of\ $A$\ . One of the many equivalent characterizations of a
weak Hilbert space given by Pisier \cite{P} we use as the definition.

\begin{definition}
A Banach space\ $X$\ is called a weak Hilbert space, if
there exist a\ $\delta > 0$\ and a\ $C \ge 1$\ such that for every finite
dimensional subspace\ $E \subseteq X$,\  there exist a
subspace\ $F \subseteq E$\ and a
projection\ $P: X \to F$\ with\ $\dim F \ge \delta \dim E$,\ $d(F) \le C$\ and
\ $||P|| \le C$\ . 
\end{definition}

The result below, which is a special case of one of the main theorems
of \cite{NTJ}, Theorem 3.11, is the main tool of this note.

\begin{theorem}
Let $X$ be a Banach lattice, which is a weak Hilbert space.
There exists a constant $K$ with the following property:\\
For every 1-unconditional, normalized sequence\ $(x_j)_{j \in I} \subseteq
X$\ (finite or infinite), and every\ $m \in \Bbb N$\ there exist a\ $J
\subseteq I$,\  $|J|=m$,\  so that if\ $A \subseteq I \backslash J$\
with\ $|A| \le m2^m$,\  then\ $(x_j)_{j \in A}$\ is K-equivalent to the unit
vector based of\ $\ell_2^{|A|}$.
\end{theorem}

We refer to \cite{LT2} for the definition of a K\"othe function space and we
let $\lambda$ denote the Lebesgue measure on [0,1]. If\ $A \subseteq
[0,1]$\ we let\ $1_A$\ denote the indicator function for $A$, and if $f$
is a real-valued function on [0,1] then we put\ $P_Af = f 1_A$\ .

Finally we refer to \cite{NTJ} for any unexplained notion on Banach lattices,
which are weak Hilbert spaces.
\vspace{1cm}

\section{The main result}

Our main result is the following

\begin{theorem}
If $X$ is a K\"othe function space on [0,1], which is a weak
Hilbert space, then $X$ is lattice isomorphic to\  $L_2(0,1)$\ .
\end{theorem}

{\bproof}
Let $K$ be the constant from Theorem 2. We wish to show that there is a 
constant $C \ge 1$ so that every finite
sequence of positive, normalized and mutually disjoint elements of $X$
is C-equivalent to the unit vector basis of a Hilbert space. It then
follows from \cite{LT2}, Theorem 1.b.13, that $X$ is lattice isomorphic to\ 
$L_2(0,1)$\ .

Hence let\ $n \in \Bbb N$\ and\ $(f_j)_{j=1}^n \subseteq X$\ be positive,
normalized with mutually disjoint supports and put\ $A_j =\mbox{supp} f_j$\ for
all\ $1 \le j \le n$\ .

Let\ $0 < \varepsilon < \frac{1}{2}$\ be arbitrary. Since $X$ is a weak
Hilbert space it is order continuous and therefore there exists a\
$\delta > 0$\ so that if\ $A \subseteq [0,1]$\ is Lebesgue measurable then

\begin{equation}
\lambda (A) \le \delta \Longrightarrow ||P_Af_j|| \le
\frac{\varepsilon}{n}\mbox{\quad \quad for all $1 \le j \le n$}.
\end{equation}

Let now $m \in \Bbb N$ be chosen, so that $m2^{-m} \le \delta$ \quad and $m
\ge n$ and put 

\begin{equation}
B_k= [(k-1)2^{-m},k2^{-m}[\mbox {\quad \quad for all $1 \le k \le 2^m$}.
\end{equation}

\begin{equation}
I=\{ (j,k)|P_kf_j \ne 0\}, \quad \quad \mbox{ where $P_k=P_{B_k},
\quad \quad 1 \le k \le n$}. 
\end{equation}

\begin{equation}
g_{jk} = ||P_kf_j||^{-1}P_kf_j\mbox{\quad \quad for all $(j,k) \in I$}.
\end{equation}

Since $g_{jk}$ is a normalized 1-unconditional sequence and\ $|I|\le
n2^m \le m2^m$\ there is a set\linebreak 
\ $J \subseteq I$\ with\ $|I \backslash J| = m$\ 
so that\ $(g_{jk})_{(j,k) \varepsilon J}$\ is $K$-equivalent to the unit
vector bases of\ $\ell_2^{|J|}$.

Put

\begin{equation}
A= \bigcup\ \{ A_j \cap B_k | (j,k) \notin J\} ,\quad \quad B=[0,1] \backslash
A.
\end{equation}

Clearly $\lambda (A) \le m2^{-m} \le \delta$ and hence

\begin{equation}
||P_Af_j|| \le \frac{\varepsilon}{n}\mbox{\quad \quad for all $1 \le j \le n$}.
\end{equation}

If $J_j =\{k|(j,k) \in J\}$ \quad \quad then

\begin{equation}
P_Bf_j = \sum_{k \in J_j}P_kf_j\quad \mbox{for all\ $1 \le j \le n$}.
\end{equation}

This shows that\ $(P_Bf_j)_{j=1}^n$\ is a block basis of\ $(g_{jk})$\ and
since\ $||P_Bf_j|| \ge 1- \varepsilon \ge \frac {1}{2}$\ for all\ $1 \le j
\le n$\ it is $2K$-equivalent to the unit vector basis of\ $\ell_2^n$\ .

Furthermore, since

\begin{equation}
2 \sum_{j=1}^n||P_Bf_j - f_j|| \le 2 \varepsilon < 1
\end{equation}

it follows that\ $(f_j)$\ is\ $2K(1+2 \varepsilon)(1-2
\varepsilon)^{-1}$-equivalent to the unit vector basis of $\ell_2^n$.
{\eproof}

As a corollary we obtain

\begin{theorem}
Let $X$ be a separable Banach lattice, which is not purely atomic and
which is not isomorphic to a Hilbert space. If $X$ is a weak Hilbert
space then there is a infinite sequence $(x_j)$ of mutually disjoint
positive normalized vectors so that $X$ is lattice isomorphic to\ 
$L_2(0,1) \oplus [x_j]$. 
\end{theorem}

\bproof
Since every separable non-atomic Banach lattice is lattice
isomorphic to a K\"othe function space on $[0,1]$ it follows from Theorem 3
and our assumptions that $X$ has atoms.
Let\ $(x_n) \subseteq X$\ be the sequence of atoms. We can then write

\begin{equation}
X = [x_j]^\perp \oplus [x_j].
\end{equation}
$[x_j]^\perp$ is a non-atomic lattice and therefore lattice isomorphic to
$L_2(0,1)$ by Theorem 3.
\eproof

It is readily verified that the argument in the proof of Theorem
3 still works, if we just assume that the space $X$ there has the
following property: There is a constant\ $K \ge 1$\ and a function\ $k:
\Bbb N \longrightarrow \Bbb N$\ with\ $k(m)m^{-1} \longrightarrow
\infty$\ for\ $m \longrightarrow \infty$\ so that every normalized
sequence in $X$
consisting of mutually disjoint elements has property\ $E(m,k(m),K)$\ for
all $m$ as defined in \cite{NTJ}.

It seems unknown whether a general Banach lattice with this property is
a weak Hilbert space, if\ $k(m)2^{-m} \longrightarrow 0$\ for\ $m
\longrightarrow \infty$\ . Compare with the proof of Theorem 4.1 in \cite{NTJ}.

\vspace{3cm}

Department of Mathematics, Odense University,
Campusvej 55, DK-5230 Odense M, Denmark.
email: njn@imada.ou.dk

\end{document}